\documentclass{amsart}
\usepackage{amssymb,color,hhline}
\usepackage{amsmath}
\usepackage{a4}
\newtheorem{theorem}{Theorem}[section]
\newtheorem{lemma}[theorem]{Lemma}
\newtheorem{remark}[theorem]{Remark}

\makeatletter
  
  \@addtoreset{equation}{section}
 \makeatother
\usepackage{color}

\def\qedbox{\hbox{$\rlap{$\sqcap$}\sqcup$}}

\begin{document}
  \title{Geometric Realizations of Hermitian curvature models}
\author{M.  Brozos-V\'azquez, P. Gilkey, H. Kang, and S. Nik\v cevi\'c}
\address{MB: Department of Mathematics, University of A Coru\~na, Spain\\
E-mail: mbrozos@udc.es}
\address{PG: Mathematics Department, University of Oregon\\
   Eugene OR 97403 USA\\
   E-mail: gilkey@uoregon.edu}
\address{HK: School of Mathematics, University of Birmingham,
Birmingham {B15 2TT, UK}\\
E-mail: {kangh@bham.ac.uk}}
\address{SN: Mathematical Institute, Sanu,
Knez Mihailova 35, p.p. 367\\
11001 Belgrade,
Serbia\\
E-mail: stanan@mi.sanu.ac.rs}
\begin{abstract} We show that a Hermitian algebraic curvature model satisfies the Gray identity if and only if it is geometrically realizable by a Hermitian
manifold. Furthermore, such a curvature model can in fact be realized by a Hermitian manifold of constant scalar curvature and constant
$\star$-scalar curvature which satisfies the Kaehler condition at the point in question.
\end{abstract}
\maketitle

\section{Introduction}

A central area of study in Differential Geometry is the examination of the relationship between purely
algebraic properties of the Riemann curvature tensor and the underlying geometric properties of the manifold.
Many authors have worked in this area in recent years. Nevertheless, many fundamental questions remain
unanswered. In this paper we shall work in the context of Hermitian geometry and refer to a few earlier works
in the almost-Hermitian context \cite{F94,MC05,MC06,Sa03} and in the Hermitian context
\cite{AGI97, B90,JK07,Sa04,V07}; the field is a vast one. Our main result, see Theorem
\ref{thm-1.4} below, gives a complete answer to the fundamental question of when a curvature tensor in a
Hermitian vector space is geometrically realizable by a Hermitian manifold.

If the almost complex structure of an almost Hermitian manifold is integrable (i.e. if the manifold is
Hermitian), then Gray
\cite{gray} showed that the curvature tensor has an additional symmetry given in Equation (\ref{eqn-1.d})
below; we refer to Theorem \ref{thm-1.3} for details. By contrast, there is no additional condition imposed
on the curvature tensor of an almost Hermitian manifold, see Theorem \ref{thm-1.1} below. Thus it is quite
striking that a geometric integrability condition imposes an additional algebraic symmetry on the curvature
tensor. In this paper, we show that every algebraic curvature tensor satisfying the Gray condition is
geometrically realized by some point $P$ of a Hermitian manifold. This result can be regarded as a converse
to Gray's result in that it shows that all the universal symmetries which the curvature tensor of Hermitian
manifold has are generated by the Gray condition and the usual curvature symmetries (see Equation
(\ref{eqn-1.a})). This result emphasizes the real difference between almost Hermitian and Hermitian
manifolds.

Let $\Omega$ be the Kaehler form. We will show, see Remark \ref{rmk-1.5} below, that the Hermitian geometric
realization can be chosen so that $d\Omega(P)=0$. Thus imposing the Kaehler identity $d\Omega (P)=0$ at a
single point imposes no additional curvature restrictions. If $d\Omega=0$ globally, then the manifold is said
to be almost Kaehler. This is a very rigid structure, see for example the discussion in
\cite{T06}, and there are additional curvature restrictions. Thus our result also emphasizes the difference
between
$d\Omega$ vanishing at a single point and
$d\Omega$ vanishing globally. Finally, we will show that the geometric realization can be chosen to have
constant scalar curvature and constant $\star$-scalar curvature. Thus, again, these conditions give rise to
no additional curvature identities.

Tricerri and Vanhecke \cite{TV81} gave a complete decomposition of the space of algebraic curvature tensors
on a Hermitian vector space into irreducible subspaces under the action of the unitary group which has been
used by many authors \cite{Bu07,F05,GM08,GM08a,MO08}. Among these factors $\mathcal{W}_7$, see Theorem
\ref{thm-2.1} for a precise definition, plays a central role in our discussion, since the space of vectors
satisfying the Gray condition is exactly
$\mathcal{W}_7^\perp$ as will be discussed in Section \ref{sect-3}. We note that either the complex Jacobi
operator or the complex curvature operator completely determine the components in
$\mathcal{W}_7^\perp$ of a curvature tensor \cite{BVGGR08}. Also note that the algebraic condition
determining $\mathcal{W}_7$ plays a role in the study of Jacobi--Ricci commuting curvature tensors \cite{GN}.

\subsection{Complex curvature models} We begin by putting things in an algebraic context. Let $V$ be a real
vector space of dimension $2n$. We shall always assume $2n\ge4$ as the
$2$-dimensional setting is trivial. Fix a symmetric positive definite inner product
$\langle\cdot,\cdot\rangle$ on $V$. Let $J$ be a {\it Hermitian complex structure} on $V$, i.e. a linear transformation of $V$
which satisfies
$$J^2=-\operatorname{id}\quad\text{and}\quad J^*\langle\cdot,\cdot\rangle=\langle\cdot,\cdot\rangle\,.$$
Let $\mathcal{U}$ be the associated unitary group;
$$\mathcal{U}:=\{U\in\operatorname{GL}(V):UJ=JU\quad\text{and}\quad U^*\langle\cdot,\cdot\rangle=\langle\cdot,\cdot\rangle\}\,.$$
Let $\mathfrak{A}\subset\otimes^4V^*$ be the space of algebraic curvature tensors; $A\in\mathfrak{A}$ if and only if $A$ satisfies the
symmetries of the Riemann curvature tensor:
\begin{equation}\label{eqn-1.a}
\begin{array}{l}
A(x,y,z,w)=-A(y,x,z,w)=A(z,w,x,y),\\
A(x,y,z,w)+A(y,z,x,w)+A(z,x,y,w)=0\,.\vphantom{\vrule height 11pt}
\end{array}\end{equation}
Let $A\in\mathfrak{A}$. We shall say that
$\mathfrak{C}:=(V,\langle\cdot,\cdot\rangle,J,A)$ is a {\it complex curvature model}.

\subsection{Almost Hermitian geometry} We say that $\mathcal{M}:=(M,g,\mathcal{J})$ is an {\it almost Hermitian manifold} if $\mathcal{J}$
is an almost complex structure on the tangent bundle $TM$ with $\mathcal{J}^*g=g$. The {\it Kaehler form}
$\Omega\in\Lambda^2(T^*M)$ is defined by $$\Omega(x,y):=g(x,\mathcal{J}y)\,.$$
Let $R$ be the Riemann curvature tensor. Let
$(T_PM,g_P,\mathcal{J}_P,R_P)$ be the associated complex curvature model for $P\in M$. A complex curvature model
$\mathfrak{C}=(V,\langle\cdot,\cdot\rangle,J,A)$ is said to be {\it geometrically realized} by $\mathcal{M}$ at $P$ if there is an isomorphism $\phi$
from $V$ to $T_PM$ so that $\langle\cdot,\cdot\rangle=\phi^*g_P$, so that $J=\phi^*\mathcal{J}_P$, and so that $A=\phi^*R_P$. We have the
following geometrical realization result in this setting
\cite{GBKNW08}:

\begin{theorem}\label{thm-1.1}
Let $\mathfrak{C}$ be a complex curvature model. There exists an almost Hermitian manifold $\mathcal{M}$
and a point $P$ so that $\mathfrak{C}$ is geometrically realized by
$\mathcal{M}$ at $P$.
\end{theorem}

\begin{remark}\label{rmk-1.2}\rm In fact, more is true. The manifold $\mathcal{M}$ in Theorem \ref{thm-1.1} can be chosen to have constant scalar
curvature and constant $\star$-scalar curvature \cite{GBKNW08}.\end{remark}

\subsection{Hermitian geometry} We say an almost Hermitian manifold $\mathcal{M}=(M,g,\mathcal{J})$ is {\it Hermitian} if $\mathcal{J}$ is an
integrable almost complex structure, i.e. the Nijenhuis tensor vanishes or, equivalently, in a neighborhood of any point of the manifold there are local
coordinates
$(x_1,...,x_n,y_1,...,y_n)$ so that
$$\mathcal{J}\partial_{x_i}=\partial_{y_i}\quad\text{and}\quad \mathcal{J}\partial_{y_i}=-\partial_{x_i}\,.$$

We introduce the following linear $\mathcal{U}$ invariant subspace of $\mathfrak{A}$:
\begin{eqnarray}
&&\mathcal{W}_G:=\{A\in\mathfrak{A}:A(x,y,z,w)+A(Jx,Jy,Jz,Jw)\nonumber\\
&&\qquad\qquad\qquad=A(Jx,Jy,z,w)+A(x,y,Jz,Jw)+A(Jx,y,Jz,w)\label{eqn-1.d}\\
&&\qquad\qquad\qquad+A(x,Jy,z,Jw)+A(Jx,y,z,Jw)+A(x,Jy,Jz,w)\}\nonumber
\end{eqnarray}
for all vectors $x,y,z,w$ in $V$. Gray \cite{gray} showed:

\begin{theorem}\label{thm-1.3}
If $P$ is a point of a Hermitian manifold $\mathcal{M}$, then $R_P\in\mathcal{W}_G$.
\end{theorem}

The following is the main result of this paper which generalizes Theorem \ref{thm-1.1} to the Hermitian
setting:

\begin{theorem}\label{thm-1.4}
Let $\mathfrak{C}=(V,\langle\cdot,\cdot\rangle,J,A)$ be a complex curvature model. Then
$A$ belongs to $\mathcal{W}_G$ if and only if there exists a Hermitian manifold $\mathcal{M}$
and a point $P$ so that $\mathfrak{C}$ is geometrically realized by
$\mathcal{M}$ at $P$.
\end{theorem}

\begin{remark}\label{rmk-1.5}
\rm We shall in fact show a bit more later. It will follow from our analysis that the Hermitian manifold $\mathcal{M}$ in Theorem
\ref{thm-1.4} can be chosen to have constant scalar curvature, to have constant $\star$-scalar curvature,
and to have $d\Omega(P)=0$ at the point in question. Thus, in
particular, imposing the Kaehler condition (i.e. $d\Omega=0$) at a single point yields no additional
curvature restrictions.
\end{remark}

\subsection{Outline of the paper}
In Section \ref{sect-2}, we present the Tricerri-Vanhecke decomposition and we linearize the problem. In Section
\ref{sect-3}, we prove Theorem \ref{thm-1.4} and verify Remark \ref{rmk-1.5}. The proof of Theorem
\ref{thm-1.1} in
\cite{GBKNW08} followed very different lines. We first realized the associated real curvature model $(V,\langle\cdot,\cdot\rangle,A)$ geometrically and then
constructed an extension of $J$ to the tangent bundle. This construction did not preserve the integrability condition on the almost complex structure and thus
very different methods are required to prove Theorem \ref{thm-1.4}. Our constructions throughout this paper will all be local; the manifolds we
will construct are neither compact nor complete. The complex analogue of the Yamabe problem for almost Hermitian manifolds has been considered in \cite{RS03}.

\section{Algebraic preliminaries}\label{sect-2}

\subsection{Decomposing $V^*\otimes V^*$ into irreducible $\mathcal{U}$ modules}
We may decompose
$$V^*\otimes V^*=S^2(V^*)\oplus\Lambda^2(V^*)$$
as the direct sum of the symmetric and of the alternating $2$-tensors, respectively.
We may use $J$ to further decompose
$$S^2(V^*)=S_+^2(V^*)\oplus S^2_-(V^*)\quad\text{and}\quad\Lambda^2(V^*)=\Lambda^2_+(V^*)\oplus\Lambda^2_-(V^*)$$
into the $\pm1$ eigenspaces of the action of $J^*$. We have a further decomposition:
$$
S_+^2(V^*)=\langle\cdot,\cdot\rangle\cdot\mathbb{R}\oplus S_{0,+}^2(V^*)\quad\text{and}\quad
\Lambda_+^2(V^*)=\Omega\cdot\mathbb{R}\oplus\Lambda_{0,+}^2(V^*)\,;
$$
$S_{0,+}^2(V^*)$ is the space of trace free symmetric $2$-tensors invariant
under $J$. One thus has a decomposition of $V^*\otimes V^*$ as the direct sum of 6 irreducible $\mathcal{U}$ modules:
\begin{equation}\label{eqn-2.a}
V^*\otimes V^*=\langle\cdot,\cdot\rangle\cdot\mathbb{R}\oplus S_{0,+}^2(V^*)\oplus S_-^2(V^*)
  \oplus\Omega\cdot\mathbb{R}\oplus\Lambda_{0,+}^2(V^*)\oplus\Lambda_-^2(V^*)\,.
\end{equation}
If $\theta$ is a
$J$-invariant symmetric form, we can define $\omega_\theta(x,y):=\theta(x,Jy)$. This defines a $\mathcal{U}$ equivariant isomorphism
$S_{0,+}^2(V^*)\approx\Lambda_{0,+}^2(V^*)$. Thus this representation and the trivial representation both appear with multiplicity $2$ in $V^*\otimes V^*$; the
other $2$ summands appear with multiplicity $1$. We have:
$$\begin{array}{ll}
\dim\{S_{0,+}^2(V^*)\}=n^2-1,&\dim\{S_-^2(V^*)\}=n^2+n,\\
\dim(\Lambda_{0,+}^2(V^*)\}=n^2-1,&\dim\{\Lambda_-^2(V^*)\}=n^2-n\,.\vphantom{\vrule height 11pt}
\end{array}$$
If $\theta\in V^*\otimes V^*$, let $\theta_{\pm,S}$, $\theta_{\pm,\Lambda}$, and $\theta_{0,+,S}$ be the
appropriate components of $\theta$ in the decomposition given above:
\begin{equation}\label{eqn-2.b}
\begin{array}{l}
\theta_{\pm,S}(x,y):=\textstyle\frac14\{\theta(x,y)+\theta(y,x)\pm\theta(Jx,Jy)\pm\theta(Jy,Jx)\},\\
\theta_{\pm,\Lambda}(x,y):=\textstyle\frac14\{\theta(x,y)-\theta(y,x)\pm\theta(Jx,Jy)\mp\theta(Jy,Jx)\},\vphantom{\vrule height 11pt}\\
$$\theta_{0,+,S}(x,y):=\theta_{+,S}(x,y)-\textstyle\frac1{2n}\{\operatorname{Tr}_{\langle\cdot,\cdot\rangle}\theta\}\langle x,y\rangle\,.
\vphantom{\vrule height 11pt}
\end{array}\end{equation}

\subsection{The scalar and $\star$-scalar curvature}\label{sect-2.3}
Adopt the {\it Einstein convention} and sum over repeated indices. Let $\{e_i\}$ be an orthonormal basis for $V$. We define the {\it Ricci
tensor}
$\rho$, the {\it scalar curvature}
$\tau$, the {\it $\star$-Ricci tensor} $\rho^\star$, and the {\it $\star$-scalar curvature} $\tau^\star$, respectively, by using the metric to contract
indices:
\begin{eqnarray*}
&&\rho(x,y):=A(e_i,x,y,e_i),\qquad\quad\tau:=\rho(e_i,e_i),\\
&&\rho^\star(x,y):=A(e_i,x,Jy,Je_i),\quad\tau^\star:=\rho^\star(e_i,e_i)\,.
\end{eqnarray*}
The Ricci tensor $\rho$ is always symmetric but the $\star$-Ricci tensor $\rho^\star$ need not be symmetric.

\subsection{The Tricerri-Vanhecke decomposition}
Give $\mathfrak{A}\subset\otimes^4V^*$ the induced inner product. Pullback defines a natural orthogonal action of
$\mathcal{U}$ on $\mathfrak{A}$. Tricerri and Vanhecke
\cite{TV81} gave a decomposition of $\mathfrak{A}$ as a $\mathcal{U}$ module:

\begin{theorem}\label{thm-2.1}
\ \begin{enumerate}
\item We have the following orthogonal direct sum decomposition of $\mathfrak{A}$ into irreducible $\mathcal{U}$ modules:
\begin{enumerate}
\item If $2n=4$,
$\mathfrak{A}=\mathcal{W}_1\oplus\mathcal{W}_2\oplus\mathcal{W}_3\oplus\mathcal{W}_4\oplus\mathcal{W}_7
\oplus\mathcal{W}_8\oplus\mathcal{W}_9$.
\item If $2n=6$,
$\mathfrak{A}=\mathcal{W}_1\oplus\mathcal{W}_2\oplus\mathcal{W}_3\oplus\mathcal{W}_4\oplus\mathcal{W}_5\oplus\mathcal{W}_7
\oplus\mathcal{W}_8\oplus\mathcal{W}_9\oplus\mathcal{W}_{10}$.
\item If $2n\ge8$,
$\mathfrak{A}=\mathcal{W}_1\oplus\mathcal{W}_2\oplus\mathcal{W}_3\oplus\mathcal{W}_4\oplus\mathcal{W}_5\oplus\mathcal{W}_6\oplus\mathcal{W}_7
\oplus\mathcal{W}_8\oplus\mathcal{W}_9\oplus\mathcal{W}_{10}$.
\end{enumerate}
We have $\mathcal{W}_1\approx\mathcal{W}_4$ and, if $2n\ge6$, $\mathcal{W}_2\approx\mathcal{W}_5$. The other $\mathcal{U}$ modules
appear with multiplicity 1.
\item We have that:
\begin{enumerate}
\item  $\tau\oplus\tau^\star:\mathcal{W}_1\oplus\mathcal{W}_4\approx\mathbb{R}\oplus\mathbb{R}$.
\item If $2n=4$, $\rho_{0,+,S}:\mathcal{W}_2\approx S_{0,+}^2(V^*)$.
\item If $2n\ge6$, $\rho_{0,+,S}\oplus\rho_{0,+,S}^\star:\mathcal{W}_2\oplus\mathcal{W}_5\approx S_{0,+}^2(V^*)\oplus S_{0,+}^2(V^*)$.
\item $\mathcal{W}_3
=\{A\in\mathfrak{A}:A(x,y,z,w)=A(Jx,Jy,z,w)\ \forall x,y,z,w\}\cap\ker(\rho)$.
\item If $2n\ge8$, $\mathcal{W}_6=\ker(\rho\oplus\rho^\star)\cap\{A\in\mathfrak{A}:J^*A=A\}\cap W_3^\perp$.
\item  $\mathcal{W}_7=\{A\in\mathfrak{A}:A(Jx,y,z,w)=A(x,y,Jz,w)\forall x,y,z,w\}$.
\item $\rho_{-,S}:\mathcal{W}_8\approx S_-^2(V^*)$.
\item  $\rho_{-,\Lambda}^*:\mathcal{W}_9\approx\Lambda_-^2(V^*)$.
\item If $2n\ge6$, $\mathcal{W}_{10}=\{A\in\mathfrak{A}:J^*A=-A\}\cap\ker(\rho\oplus\rho^\star)$.
\end{enumerate}
\def\gronk{\noalign{\hrule}\vphantom{\vrule height 13pt\frac1{4_{A_A}}}}
\item The dimensions of these modules are given by:\smallbreak\noindent\qquad{$\begin{array}{|l|l|l|l|l|}\gronk
&\dim(V)=4&\dim(V)=6&\dim(V)=2n\ge8\\\gronk
\mathcal{W}_1&1&1&1\\\gronk
\mathcal{W}_2&3&8&n^2-1\\\gronk
\mathcal{W}_3&5&27&\frac14n^2(n-1)(n+3)\\\gronk
\mathcal{W}_4&1&1&1\\\gronk
\mathcal{W}_5&0&8&n^2-1\\\gronk
\mathcal{W}_6&0&0&\frac14n^2(n+1)(n-3)\\\gronk
\mathcal{W}_7&2&12&\frac16n^2(n^2-1)\\\gronk
\mathcal{W}_8&6&12&n^2+n\\\gronk
\mathcal{W}_9&2&6&n^2-n\\\gronk
\mathcal{W}_{10}&0&30&\frac23n^2(n^2-4)\\\noalign{\hrule}
\end{array}$}
\end{enumerate}
\end{theorem}

The following will be a useful observation in what follows.

\begin{remark}\label{rmk-2.2}
\rm Let $\mathfrak{B}$ be a $\mathcal{U}$ invariant
subspace of $\mathfrak{A}$. Let $i\in\{3,6,7,8,9,10\}$. Then $\mathcal{W}_i$ appears with multiplicity $1$ in $\mathfrak{A}$ as a $\mathcal{U}$ module.
Consequently we either have that
$\mathcal{W}_i\subset\mathfrak{B}$ or that $\mathfrak{B}\subset\mathcal{W}_i^\perp=\oplus_{j\ne i}
\mathcal{W}_j$. Let $\pi_i$ denote orthogonal projection on
$\mathcal{W}_i$. If $\pi_i\mathfrak{B}\ne\{0\}$, then $\mathcal{W}_i\subset\mathfrak{B}$.
\end{remark}

\subsection{Linearizing the problem}\label{sect-2.4}
Let $\Theta\in S_+^2(V^*)\otimes S^2(V^*)$. We define:
$$\mathcal{L}(\Theta)(x,y,z,w):=\Theta(x,z,y,w)+\Theta(y,w,x,z)-\Theta(x,w,y,z)-\Theta(y,z,x,w)\,.$$
It is easily verified that the relations of Equation (\ref{eqn-1.a}) are satisfied and thus $\mathcal{L}$ defines a $\mathcal{U}$ equivariant map
$$\mathcal{L}:S_+^2(V^*)\otimes S^2(V^*)\rightarrow\mathfrak{A}\,.$$
We let
$\mathfrak{L}=\mathfrak{L}(V,\langle\cdot,\cdot\rangle,J):=\operatorname{Range}(\mathcal{L})$ be a $\mathcal{U}$ invariant subspace.

\begin{lemma}\label{lem-2.3}
\ \begin{enumerate}
\item Let $\mathcal{M}=(M,g,J)$ be a Hermitian manifold. Let $P\in M$. Then $d\Omega(P)=0$ if and only if there exist
holomorphic coordinates
$(z_1,...,z_n)$ for $M$ centered at $P$ so that
$g(P)=\delta+O(|z|^2)$.
\item $A\in\mathfrak{L}$ if and only if the complex model $\mathfrak{C}=(V,\langle\cdot,\cdot\rangle,J,A)$ can be geometrically realized by a point $P$ of a
Hermitian manifold
$\mathcal{M}$ where $d\Omega(P)=0$.
\end{enumerate}
\end{lemma}

\begin{proof} Although Assertion (1) is well known, see for example \cite{G73}, we present the proof as this result is central to our discussion and to keep this
paper as self-contained as possible.
 Let $T_{\mathbb{C}}(M):=T(M)\otimes_{\mathbb{R}}\mathbb{C}$ be the complex tangent bundle of a Hermitian manifold $\mathcal{M}$. Extend the
metric $g$ to
$T_{\mathbb{C}}(M)$ to be complex linear in the first argument and conjugate linear in the
second argument. Choose local holomorphic coordinates $z^a=x^a+\sqrt{-1}y^a$. We have
$$\partial_{z_a}:=\textstyle\frac12(\partial_{x_a}-\sqrt{-1}\partial_{y_a})\quad\text{and}\quad
  \partial_{\bar z_a}:=\textstyle\frac12(\partial_{x_a}+\sqrt{-1}\partial_{y_a})\,.$$
The condition that $J$ is compatible with $g$ then means $g(\partial_{z_a},\partial_{\bar z_b})=0$. We set
$$
g_{a\bar b}:=g(\partial_{z_a},\partial_{z_b})=
\textstyle\frac12\{g(\partial_{x_a},\partial_{x_b})-\sqrt{-1}g(\partial_{x_a},\partial_{y_b})\}
=\bar g_{b\bar a}\,.$$
The Kaehler form is given by $\Omega=\frac{\sqrt{-1}}2g_{b\bar d}dz^b\wedge dz^{\bar d}$. Consequently:
\begin{eqnarray*}
d\Omega&=&{\textstyle\frac{\sqrt{-1}}2}\sum_{b<c,d}(g_{c\bar d/b}-g_{b\bar d/c})dz^b\wedge dz^c\wedge d\bar z^d\\
&-&
{\textstyle\frac{\sqrt{-1}}2}
\sum_{b,c<d}(g_{b\bar d/\bar c}-g_{b\bar c/\bar d})dz^b\wedge d\bar z^c\wedge d\bar z^d\,.
\end{eqnarray*}
This shows that the condition $d\Omega(P)=0$ is equivalent to the symmetry:
\begin{equation}\label{eqn-2.c}
g_{b\bar d/c}(P)=g_{c\bar d/b}(P)\,.
\end{equation}
Clearly if we can choose holomorphic coordinates so all the $1$-jets of the metric vanish at $P$, then Equation
(\ref{eqn-2.c}) is satisfied so $d\Omega(P)=0$. Conversely, suppose Equation (\ref{eqn-2.c}) is satisfied. Choose coordinates with $P=0$ so that
$g_{a\bar b}(0)=\delta_{a\bar b}$. Consider the holomorphic change of coordinates:
$$z^a:=w^a+\xi_{abc}w^bw^c\quad\text{where}\quad\xi_{abc}=\xi_{acb}\in\mathbb{C}\,.$$
At the point $P=0$, we may then express
\begin{eqnarray*}
&&\partial_{w_c}=\textstyle\frac{\partial z^a}{\partial w_c}\partial_{z_a}
   =\partial_{z_c}+2\xi_{abc}w^b\partial_{z_a},\\
&&g^w_{c\bar d}=g^z_{c\bar d}+2\xi_{dbc}w^b+O(\bar w)+O(|w|^2),\\
&&g^w_{c\bar d/b}=g^z_{c\bar d/b}+2\xi_{dbc}+O(|w|)\,.
\end{eqnarray*}
We set $\xi_{dbc}:=-\frac12g^z_{c\bar d/b}$; this is symmetric in $\{b,c\}$ by Equation (\ref{eqn-2.c}) and
thus defines an admissible change of coordinates so that $\partial_{w_b}g^w_{c\bar d}(0)=0$. Taking the complex
conjugate yields $\partial_{\bar w_b}g^w_{d\bar c}(0)=0$ as well and thus $dg^w(0)=0$; this establishes Assertion (1).

To prove Assertion (2), we argue as follows.
Suppose first
$A\in\mathfrak{L}$. Give
$\mathbb{R}^{2n}=\mathbb{C}^n$ the usual system of coordinates $(u^1,...,u^{2n})=(x^1,...,x^n,y^1,...,y^n)$ and integrable complex structure
$\mathcal{J}(\partial_{x_i})=\partial_{y_i}$ and $\mathcal{J}(\partial_{y_i})=-\partial_{x_i}$. Choose $\Theta\in S_+^2(V^*)\otimes S^2(V^*)$ with
$\mathcal{L}(\Theta)=A$ and define
\begin{equation}\label{eqn-2.d}
g_{ij}:=\delta_{ij}+2\Theta_{ijkl}u^ku^l\,.
\end{equation}
Since $\Theta(x,y,z,w)=\Theta(Jx,Jy,z,w)$, $\mathcal{J}^*g=g$. Let $B_\epsilon$ be the Euclidean ball of radius $\epsilon>0$ centered at the origin. Since $g$ is
non-singular at the origin, there exists $\varepsilon>0$ so $g$ is non singular on $B_\varepsilon$; let $\mathcal{M}:=(B_\epsilon,g,\mathcal{J})$ be the
resulting Hermitian manifold. Since the first derivatives of the metric vanish at the origin, $d\Omega(P)=0$ and we may compute:
\begin{eqnarray*}
R(\partial_{u_i},\partial_{u_j},\partial_{u_k},\partial_{u_l})
&=&\textstyle\frac12\{\partial_{u_i}\partial_{u_k}g_{jl}+\partial_{u_j}\partial_{u_l}g_{ik}
-\partial_{u_i}\partial_{u_l}g_{jk}-\partial_{u_j}\partial_{u_k}g_{il}\}\\
&=&\Theta_{ikjl}+\Theta_{jlik}-\Theta_{iljk}-\Theta_{jkil}=A\,.
\end{eqnarray*}
This shows every element of $\mathfrak{L}$ can be represented by a Hermitian manifold with $d\Omega(P)=0$. Conversely, suppose given a Hermitian
manifold
with $d\Omega(P)=0$. Choose a holomorphic coordinate system centered at $P$ where all the first derivatives of the metric vanish and so the coordinate frame is
orthonormal at $P$. Thus up to second order, $g$ has the form given in Equation (\ref{eqn-2.d}). The above calculation then shows
$R_P=\mathcal{L}(\Theta)\in\mathfrak{L}$.\end{proof}

\section{The proof of Theorem \ref{thm-1.4}}\label{sect-3}

The results of Section \ref{sect-2} reduce the proof of Theorem \ref{thm-1.4} to the assertion:
 $$\mathfrak{L}=\mathcal{W}_G\,.$$
We begin our study with the following result:

\begin{lemma}\label{lem-3.1}
 $\mathfrak{L}\subset\mathcal{W}_G\subset\mathcal{W}_7^\perp$.
\end{lemma}

\begin{proof} By Lemma \ref{lem-2.3}, every element of $\mathcal{L}$ can be geometrically realized by a Hermitian manifold. Theorem \ref{thm-1.3} now implies
$\mathfrak{L}\subset\mathcal{W}_G$.  By Remark \ref{rmk-2.2}, we may show $\mathcal{W}_G\subset\mathcal{W}_7^\perp$ by showing
$\mathcal{W}_G\cap\mathcal{W}_7=\{0\}$. Let $A\in\mathcal{W}_G\cap\mathcal{W}_7$. Since $A\in\mathcal{W}_7$, the curvature symmetries imply additionally that
\medbreak\qquad
$A(Jx,y,z,w)=-A(Jx,y,w,z)=-A(x,y,Jw,z)=A(x,y,z,Jw)$
\medbreak\qquad\qquad$=-A(y,x,z,Jw)=-A(Jy,x,z,w)=A(x,Jy,z,w)$.
\medbreak\noindent
Since $A\in\mathcal{W}_G$, Equation (\ref{eqn-1.d}) implies
\begin{eqnarray*}
&&2A(x,y,z,w)=A(x,y,z,w)+A(JJx,y,z,JJw)\\
&&\qquad=A(x,y,z,w)+A(Jx,Jy,Jz,Jw)\\
&&\qquad=A(Jx,Jy,z,w)+A(x,y,Jz,Jw)+A(Jx,y,Jz,w)\\
&&\qquad+A(x,Jy,z,Jw)+A(Jx,y,z,Jw)+A(x,Jy,Jz,w)\\
&&\qquad=-6A(x,y,z,w)\,.
\end{eqnarray*}
Consequently $8A(x,y,z,w)=0$ so $A=0$.
\end{proof}

We continue our study with:

\begin{lemma}\label{lem-3.2}
\ \begin{enumerate}
\item $\tau\oplus\tau^\star:\mathfrak{L}\rightarrow\mathbb{R}\oplus\mathbb{R}\rightarrow0$. Thus
$\mathcal{W}_1\oplus\mathcal{W}_4\subset\mathfrak{L}$.
\item If $2n=4$, then $\rho_{0,+,S}^\star:\mathfrak{L}\rightarrow S_{0,+}^2(V^*)\rightarrow0$. Thus $\mathcal{W}_2\subset\mathfrak{L}$.
\item $\rho_{-,S}:\mathfrak{L}\rightarrow S_-^2(V^*)\rightarrow0$. Thus $\mathcal{W}_8\subset\mathfrak{L}$.
\item $\rho_{-,\Lambda}^*:\mathfrak{L}\rightarrow\Lambda^2_-(V^*)\rightarrow0$. Thus
$\mathcal{W}_9\subset\mathfrak{L}$.
\item If $2n\ge6$, then $\{\rho_{0,+,S}\oplus\rho_{0,+,S}^\star\}:\mathfrak{L}\rightarrow\{S_{0,+}^2(V^*)\oplus S_{0,+}^2(V^*)\}\rightarrow0$. Thus\newline
$\mathcal{W}_2\oplus\mathcal{W}_5\subset\mathfrak{L}$.
\item $\mathfrak{L}\cap\mathcal{W}_3\ne\{0\}$. Thus $\mathcal{W}_3\subset\mathfrak{L}$.
\item $\mathfrak{L}\cap\mathcal{W}_{10}\ne\{0\}$. Thus $\mathcal{W}_{10}\subset\mathfrak{L}$.
\item If $2n\ge6$, then $\mathfrak{L}\cap\mathcal{W}_6\ne\{0\}$. Thus $\mathcal{W}_6\subset\mathfrak{L}$.
\end{enumerate}\end{lemma}

\begin{proof} We shall use Remark \ref{rmk-2.2} to examine the relationship of $\mathcal{W}_i$ to $\mathfrak{L}$ for the indices $i\in\{3,6,7,8,9,10\}$; the
subspaces
$\mathcal{W}_1$ and
$\mathcal{W}_2$ can appear with multiplicity $2$ and thus require slightly more care. We use the formalism of Lemma \ref{lem-2.3}. We shall define
metrics $g=\delta+O(|x|^2)$ and let $A\in\mathfrak{L}$ be the curvature tensor at the origin. Let
$$\xi\circ\eta:=\textstyle\frac12(\xi\otimes\eta+\eta\otimes\xi)$$
denote the symmetric product. Let
$\varrho$ and
$\varepsilon$ be real constants. Consider the Hermitian metric:
$$g=\delta-\varepsilon x_1^2(dx_1\circ dx_1+dy_1\circ dy_1)-\varrho x_1^2(dx_2\circ dx_2+dy_2\circ dy_2)\,.$$
The non-zero curvatures then become, up to the usual $\mathbb{Z}_2$ symmetries,
$$A(\partial_{x_1},\partial_{y_1},\partial_{y_1},\partial_{x_1})=\varepsilon,\quad
  A(\partial_{x_1},\partial_{x_2},\partial_{x_2},\partial_{x_1})=\varrho,\quad
  A(\partial_{x_1},\partial_{y_2},\partial_{y_2},\partial_{x_1})=\varrho\,.
$$
This shows
$\tau=2\varepsilon+4\varrho$ and $\tau^\star=2\varepsilon$ so $\tau\oplus\tau^\star$ is a surjective map to $\mathbb{R}\oplus\mathbb{R}$.
Thus Assertion (1) follows from Theorem \ref{thm-2.1}:
$$\mathcal{W}_1\oplus\mathcal{W}_4\subset\mathfrak{L}\,.$$

The non-zero entries in the Ricci tensor are given by:
$$\begin{array}{ll}
\rho(\partial_{x_1},\partial_{x_1})=\varepsilon+2\varrho,&\rho(\partial_{y_1},\partial_{y_1})=\varepsilon,\\
  \rho(\partial_{x_2},\partial_{x_2})=\varrho,&\rho(\partial_{y_2},\partial_{y_2})=\varrho\,.
\end{array}$$
We take $\varrho=-1$ and $\varepsilon=2$ to ensure $\rho$ is trace free. We use Equation (\ref{eqn-2.b}) to see
$$\begin{array}{ll}
\rho_{0,+,S}(\partial_{x_1},\partial_{x_1})=1,&\rho_{0,+,S}(\partial_{y_1},\partial_{y_1})=1,\\
\rho_{0,+,S}(\partial_{x_2},\partial_{x_2})=-1,&\rho_{0,+,S}(\partial_{y_2},\partial_{y_2})=-1\,.\vphantom{\vrule height 11pt}
\end{array}$$
This shows that $\rho_{0,+,S}$ is non-zero on $\mathfrak{L}$; Assertion (2) now follows if $2n=4$ since $\mathcal{W}_5$ is not present:
$$\mathcal{W}_2\subset\mathfrak{L}\,.$$

 We use Equation
(\ref{eqn-2.b}) to compute similarly
$$\begin{array}{ll}
\rho_{-,S}(\partial_{x_1},\partial_{x_1})=-1,&\rho_{-,S}(\partial_{y_1},\partial_{y_1})=1,\\
\rho_{-,S}(\partial_{x_2},\partial_{x_2})=0,&\rho_{-,S}(\partial_{y_2},\partial_{y_2})=0\,.\vphantom{\vrule height 11pt}
\end{array}$$
This shows $\rho_{-,S}$ is non-trivial on $\mathfrak{L}$ and Assertion (3) follows:
$$\mathcal{W}_8\subset\mathfrak{L}\,.$$

We clear the previous notation and consider:
$$g=\delta-2\varepsilon x_1^2(dx_1\circ dx_2+dy_1\circ dy_2)\,.$$
There is only one non-zero curvature entry $A(\partial_{x_1},\partial_{y_1},\partial_{y_2},\partial_{x_1})=2\varepsilon$. We let
$A^*(x,y,z,w)=A(x,y,Jz,Jw)$ and use Equation (\ref{eqn-2.b}) to compute:
\begin{eqnarray*}
&&A^*(\partial_{x_1},\partial_{y_1},\partial_{x_2},\partial_{y_1})=A^*(\partial_{y_2},\partial_{x_1},\partial_{y_1},\partial_{x_1})=-2\varepsilon,\\
&&\rho^\star(\partial_{x_1},\partial_{x_2})=\rho^\star(\partial_{y_2},\partial_{y_1})=2\varepsilon,\\
&&\rho^\star_\Lambda(\partial_{x_1},\partial_{x_2})=-\rho^\star_\Lambda(\partial_{x_2},\partial_{x_1})=\rho^\star_\Lambda(\partial_{y_2},\partial_{y_1})=
-\rho^\star_\Lambda(\partial_{y_1},\partial_{y_2})=\varepsilon\,.
\end{eqnarray*}
This shows $0\ne\rho^\star_\Lambda\in\Lambda^2_-$ so $A$ has a non-trivial component in $\mathcal{W}_9$. This
completes the proof of Assertion (4):
$$\mathcal{W}_9\subset\mathfrak{L}\,.$$

Assume $2n\ge6$. We clear the
previous notation and consider:
$$
g=\delta-2\varrho x_1^2(dx_1\circ dx_2+dy_1\circ dy_2)-2\varepsilon x_1^2(dx_2\circ dx_3+dy_2\circ dy_3)\,.
$$
The non-zero curvatures now become:
$$
A(\partial_{x_1},\partial_{y_1},\partial_{y_2},\partial_{x_1})=\varrho,\quad
A(\partial_{x_1},\partial_{x_2},\partial_{x_3},\partial_{x_1})=
A(\partial_{x_1},\partial_{y_2},\partial_{y_3},\partial_{x_1})=\varepsilon\,.
$$
Note that $\rho$ is always symmetric. We have
$$\begin{array}{ll}
\rho(\partial_{y_1},\partial_{y_2})=\varrho,&
\rho_{0,+,S}(\partial_{y_1},\partial_{y_2})=\rho_{0,+,S}(\partial_{x_1},\partial_{x_1})=\textstyle\frac12\varrho,\\
\rho(\partial_{x_2},\partial_{x_3})=\rho(\partial_{y_2},\partial_{y_3})=\varepsilon,&
\rho_{0,+,S}(\partial_{x_2},\partial_{x_3})=\rho_{0,+,S}(\partial_{y_2},\partial_{y_3})=\varepsilon\,.
\end{array}$$
We have:
$$\begin{array}{ll}
A^*(\partial_{x_1},\partial_{x_2},\partial_{y_3},\partial_{y_1})=\varepsilon,&
A^*(\partial_{x_3},\partial_{x_1},\partial_{y_1},\partial_{y_2})=\varepsilon,\\
A^*(\partial_{x_1},\partial_{y_2},\partial_{x_3},\partial_{y_1})=-\varepsilon,&
A^*(\partial_{y_3},\partial_{x_1},\partial_{y_1},\partial_{x_2})=-\varepsilon,\\
A^*(\partial_{x_1},\partial_{y_1},\partial_{x_2},\partial_{y_1})=-\varrho,&A^*(\partial_{y_2},\partial_{x_1},\partial_{y_1},\partial_{x_1})=-\varrho\,.
\end{array}$$
This shows that:
$$\begin{array}{ll}
\rho^\star(\partial_{x_1},\partial_{x_2})=\varrho,&\rho^\star(\partial_{y_2},\partial_{y_1})=\varrho,\\
\rho_{0,+,S}^\star(\partial_{x_1},\partial_{x_2})=\textstyle\frac12\varrho,&\rho_{0,+,S}^\star(\partial_{y_1},\partial_{y_2})=\textstyle\frac12\varrho\,.
\end{array}$$
If we take $\varrho=0$ and $\varepsilon\ne0$, then $\rho_{0,+,S}\ne0$ and $\rho_{0,+,S}^*=0$. Thus
\begin{eqnarray*}
&&\{S_{0,+}^2(V^*)\oplus 0\}\cap\{\rho_{0,+,S}\oplus\rho_{0,+,S}^*\}\mathfrak{L}\ne\{0\}\quad\text{so}\\
&&\{S_{0,+}^2(V^*)\oplus 0\}\subset\{\rho_{0,+,S}\oplus\rho_{0,+,S}^*\}\mathfrak{L}\,.
\end{eqnarray*}
On the other hand, if we take $\varrho\ne0$, then $\rho_{0,+,S}^*\ne0$. Thus we have a non-zero component in the second factor and
$$\{S_{0,+}^2(V^*)\oplus S_{0,+}^*(V^*)\}\subset\{\rho_{0,+,S}\oplus\rho_{0,+,S}^*\}\mathfrak{L}\,.$$
This establishes Assertion (5):
$$\mathcal{W}_2\oplus\mathcal{W}_5\subset\mathfrak{L}\,.$$

 To prove Assertion (6), we consider the metric
$$
g=\delta-2\{x_1^2+y_1^2-x_2^2-y_2^2\}(dx_1\circ dx_2+dy_1\circ dy_2)\,.
$$
The non-zero components of $A$ are then given, up to the usual $\mathbb{Z}_2$ symmetries by:
\begin{eqnarray*}
&&A(\partial_{x_1},\partial_{y_1},\partial_{y_2},\partial_{x_1})=A(\partial_{y_1},\partial_{x_1},\partial_{x_2},\partial_{y_1})=1,\\
&&A(\partial_{x_2},\partial_{y_1},\partial_{y_2},\partial_{x_2})=A(\partial_{y_2},\partial_{x_1},\partial_{x_2},\partial_{y_2})=-1,\,.
\end{eqnarray*}
We have
$\rho=0$ and $A(Jx,Jy,z,w)=A(x,y,z,w)$ for all $x$, $y$, $z$, and $w$. This shows $A\in\mathcal{W}_3$ and proves Assertion (6) by showing
$$\mathcal{W}_3\subset\mathfrak{L}\,.$$

Let $2n\ge6$. We consider
$$
g=\delta-2\{x_1^2-y_1^2\}(dx_2\circ dx_3+dy_2\circ dy_3) \,.
$$
The non-zero curvatures are then
\begin{eqnarray*}
&&A(\partial_{x_1},\partial_{x_2},\partial_{x_3},\partial_{x_1})=A(\partial_{x_1},\partial_{y_2},\partial_{y_3},\partial_{x_1})=1,\\
&&A(\partial_{y_1},\partial_{x_2},\partial_{x_3},\partial_{y_1})=A(\partial_{y_1},\partial_{y_2},\partial_{y_3},\partial_{y_1})=-1\,.
\end{eqnarray*}
This tensor has vanishing Ricci and $\star$-Ricci curvature. Since $J^*A=-A$, $A\in\mathcal{W}_{10}$. This proves Assertion (7) by showing
$$\mathcal{W}_{10}\subset\mathfrak{L}\,.$$

Let $2n\ge8$. We take
$$ds^2=\delta-2\{x_1x_2+y_1y_2\}(dx_3\circ dx_4+dy_3\circ dy_4)\,.$$
The non-zero curvatures are
\begin{eqnarray*}
&&A(\partial_{x_1},\partial_{x_3},\partial_{x_4},\partial_{x_2})=
A(\partial_{y_1},\partial_{x_3},\partial_{x_4},\partial_{y_2})\\
&=&A(\partial_{x_1},\partial_{x_4},\partial_{x_3},\partial_{x_2})=
A(\partial_{y_1},\partial_{x_4},\partial_{x_3},\partial_{y_2})\\
&=&A(\partial_{x_1},\partial_{y_3},\partial_{y_4},\partial_{x_2})=
A(\partial_{y_1},\partial_{y_3},\partial_{y_4},\partial_{y_2})\\
&=&A(\partial_{x_1},\partial_{y_4},\partial_{y_3},\partial_{x_2})=
A(\partial_{y_1},\partial_{y_4},\partial_{y_3},\partial_{y_2})=1\,.
\end{eqnarray*}
We observe that $\rho=\rho^\star=0$. Since $A(Jx,Jy,z,w)\ne A(x,y,z,w)$, $A\notin\mathcal{W}_3$. Thus $A$ has a non-zero component in
$\mathcal{W}_6\oplus\mathcal{W}_7$. As $\mathcal{L}\perp\mathcal{W}_7$, $A$ has a non-zero component in $\mathcal{W}_6$ and Assertion (8)
follows; $\mathcal{W}_6\subset\mathfrak{L}$.
\end{proof}

\medbreak\noindent{\it Proof of Theorem \ref{thm-1.4}.} By Lemma~\ref{lem-3.1}, we have
$\mathfrak{L}\subset\mathcal{W}_G\subset\mathcal{W}_7^\perp$. The assertion
$\mathcal{W}_7^\perp\subset\mathfrak{L}$ follows from the Tricerri-Vanhecke decomposition described in
Theorem \ref{thm-2.1} and from Lemma
\ref{lem-3.2}.\hfill\qedbox

\medbreak\noindent{\it Proof of Remark \ref{rmk-1.5}}. The construction given above yields $\mathcal{M}$ with $d\Omega(P)=0$ realizing the given complex
curvature model $\mathfrak{C}$ at $P$. In
\cite{GBKNW08}, we considered a further variation
$$h:=g+2\xi(dx_1\circ dx_1+dy_1\circ dy_1)+2\eta(dx_2\circ dx_2+dy_2\circ dy_2)$$
where $\{\xi,\eta\}$ are smooth functions vanishing to second order at $P$. We showed it was possible to choose $\{\xi,\eta\}$ so that the resulting metric had
constant scalar curvature and constant $\star$-scalar curvature. Since $\{\xi,\eta\}$ vanish to second order, $(M,\mathcal{J},h)$ realizes $\mathfrak{C}$ at $P$ as
well and $d\Omega_{\xi,\eta}=0$. This establishes Remark
\ref{rmk-1.5}.
\hfill\qedbox

\section*{Acknowledgments} The research of all authors partially supported by Project MTM2006-01432
(Spain). The research of P. Gilkey was also partially supported by  Project DGI SEJ2007-67810a (Spain).  Research of H. Kang was also partially supported by the
University of Birmingham (UK). Research of S. Nik\v cevi\'c was also partially supported by Research of Project 144032 (Serbia). It is also
a pleasure to acknowledge with gratitude useful conversations with G. Weingart on this subject.

\end{document}